\documentclass[a4paper,12pt]{amsart}
\usepackage{ifthen}
\usepackage{graphicx}
\usepackage{mathrsfs}
\usepackage{color}
\nonstopmode \numberwithin{equation}{section}
\newtheorem{thm}{Theorem}[section]
\newtheorem{cor}[thm]{Corollary}
\newtheorem{lem}[thm]{Lemma}

\theoremstyle{definition}

\newenvironment{pf}[1][]{%
 \vskip 3mm
 \noindent
 \ifthenelse{\equal{#1}{}}%
  {{\slshape Proof. }}%
  {{\slshape #1.} }%
 }%
{\qed\bigskip}

\newcounter{alphabet}
\newcounter{tmp}
\newenvironment{Thm}[1][]{\refstepcounter{alphabet}%
\bigskip%
\noindent%
{\bf Theorem \Alph{alphabet}}%
\ifthenelse{\equal{#1}{}}{}{ (#1)}%
{\bf .} \itshape}{\vskip 8pt}

\newcommand{\A}{{\mathcal A}}

\newcommand{\C}{{\mathbb C}}
\newcommand{\D}{{\mathbb D}}

\newcommand{\K}{{\mathcal K}}

\newcommand{\R}{{\mathbb R}}

\newcommand{\es}{{\mathcal S}}

\renewcommand{\Im}{{\,\operatorname{Im}\,}}

\renewcommand{\Re}{{\,\operatorname{Re}\,}}

\newcommand{\inv}{^{-1}}

\renewcommand{\arg}{\,{\operatorname{arg}\,}}

\newcommand{\aand}{{\quad\text{and}\quad}}

\newcounter{minutes}
\divide\time by 60
\newcounter{hours}
\multiply\time by 60 \addtocounter{minutes}{-\time}

\begin{document}
\bibliographystyle{amsplain}
\title{
Notes on convex functions of order $\alpha$}

\author[T. Sugawa]{Toshiyuki Sugawa}
\address{Graduate School of Information Sciences \\
Tohoku University \\ 
Aoba-ku, Sendai 980-8579, Japan}
\email{sugawa@math.is.tohoku.ac.jp}
\author[L.-M.~Wang]{Li-Mei Wang}
\address{School of Statistics,
University of International Business and Economics, No.~10, Huixin
Dongjie, Chaoyang District, Beijing 100029, China}
\email{wangmabel@163.com}

\begin{abstract}
Marx and Strohh\"acker showed around in 1933 that
$f(z)/z$ is subordinate to $1/(1-z)$ for a normalized convex function $f$
on the unit disk $|z|<1.$
Brickman, Hallenbeck, MacGregor and Wilken proved in 1973 further that
$f(z)/z$ is subordinate to $k_\alpha(z)/z$ if $f$ is
convex of order $\alpha$ for $1/2\le\alpha<1$
and conjectured that this is true also for $0<\alpha<1/2.$
Here, $k_\alpha$ is the standard extremal function in the class of
normalized convex functions of order $\alpha$ and $k_0(z)=z/(1-z).$
We prove the conjecture and study geometric properties
of convex functions of order $\alpha.$
In particular, we prove that $(f+g)/2$ is starlike whenever
$f$ and $g$ both are convex of order $3/5.$
\end{abstract}
\thanks{
The present research was
supported by National Natural Science Foundation of China (No.~11326080)
and JSPS Grant-in-Aid for Scientific Research (B) 22340025.} 
\maketitle

\section{Introduction and main result}
Let $\A$ denote the set of analytic functions on the open unit disk
$\D=\{z\in\C: |z|<1\}.$
Let $\A_1$ be the subclass of $\A$ consisting of functions $f$
normalized by $f(0)=f'(0)-1=0.$
Further let $\es$ be the subset of
$\A_1$ consisting of functions $f$ univalent on $\D.$
The present paper mainly deals with the subfamily of $\es,$
denoted by $\K(\alpha)$, consisting of convex functions of order $\alpha$
introduced by Robertson \cite{Rob36}.
Here, for a constant $0\le\alpha<1,$
a function $f$ in $\A_1$ is called {\it convex of order $\alpha$} if
$$
\Re\left(1+\frac{zf''(z)}{f'(z)}\right)>\alpha
$$
for $z\in \D.$
Note that the class $\K(0)=\K$ is known to consist of convex
functions in $\A_1.$
Here, a function $f$ in $\A$ is called {\it convex} if $f$
maps $\D$ univalently onto a convex domain.
A function $f\in \A$ is called starlike if $f$ maps $\D$ univalently onto
a domain starlike with respect to $f(0).$ It is clear that every
convex function is starlike.
We denote by $\es^*$ the set of starlike functions in $\A_1.$
By definition, it is obvious that for $0\leq \alpha<\beta< 1,$
$$
\K(\beta)\subset\K(\alpha)\subset\K\subset\es^*\subset\es.
$$

The Koebe function $z/(1-z)^2$ is often extremal in $\es^*$ or even in $\es$
and thus plays quite an important role in the theory of univalent functions.
It is helpful in many respects to have such an extremal function
for the class $\K(\alpha).$
Since the function $(1+(1-2\alpha)z)/(1-z)$ maps $\D$ univalently onto
the half-plane $\Re w>\alpha,$ indeed, the function
$k_\alpha\in\K(\alpha)$ characterized by the following relations serves
as an extremal one:
$$
1+\frac{zk_\alpha''(z)}{k_\alpha'(z)}=\frac{1+(1-2\alpha)z}{1-z},
\aand
k_\alpha(0)=0,~ k_\alpha'(0)=1.
$$
It is easy to find an explicit form of $k_\alpha$ as follows:
\begin{align*}
k_{\alpha}(z)
&=\begin{cases}
(1-2\alpha)\inv[(1-z)^{2\alpha-1}-1],&\quad \alpha\not=1/2,\\
-\log(1-z),&\quad \alpha=1/2.
\end{cases}
\end{align*}
We now recall the notion of subordination between
two analytic functions $f$ and $g$ on $\D.$
We say that $f$ is subordinate to $g$ and write $f\prec g$ or
$f(z)\prec g(z)$ for it if there exists an analytic
function $\omega$ on $\D$ such that $\omega(0)=0,$
$|\omega(z)|<1$ and $f(z) = g(\omega(z))$ for $z\in\D.$
When $g$ is univalent, $f$ is subordinate to $g$ precisely if $f(0) =
g(0)$ and if $f(\D)\subset g(\D).$

In 1973, Brickman, Hallenbeck, MacGregor and Wilken proved in
\cite[Theorem 11]{BHMW73} the following result for convex functions of
order $\alpha$.

\begin{Thm}[Brickman {\it et al.}]
If $f\in\K(\alpha)$ for $1/2\leq\alpha<1$, then
$$
\frac{f(z)}{z}\prec\frac{k_{\alpha}(z)}{z} \quad\text{on}~ \D.
$$
\end{Thm}

We note that $k_0(z)/z=1/(1-z)$ maps $\D$ univalently onto the half-plane 
$\Re w>1/2.$
Thus the above relation also holds when $\alpha=0$ by a theorem of
Marx and Strohh\"acker (see \cite[Theorem 10]{BHMW73}).
In \cite{BHMW73}, they conjectured that the assertion of Theorem A would hold for
$0<\alpha<1/2$ as well.
They also observed that the conjecture is confirmed if one could show that
the function $k_{\alpha}(z)/z$ is convex.
They prove the last theorem by showing it for $1/2\le\alpha<1$ 
(cf.~\cite[Lemma 3]{BHMW73}).
We will show it for all $\alpha.$

\begin{thm}\label{thm:convex}
The function $h_\alpha(z)=k_\alpha(z)/z$ maps $\D$ univalently onto a convex
domain for each $0\le\alpha<1.$
\end{thm}

We remark that, in the context of the hypergeometric function,
this follows also from results of K\"ustner in
\cite{Kus02} (see the remark at the end of Section 2 for more details).
Anyway, the conjecture has been confirmed:

\begin{cor}\label{cor:main}
Let $0\le\alpha<1.$
Then, for $f\in\K(\alpha),$ the following subordination holds:
$$
\frac{f(z)}{z}\prec\frac{k_{\alpha}(z)}{z} \quad\text{on}~\D.
$$
\end{cor}

In view of the form, it is easy to see that $k_\alpha$ is bounded on $\D$
if and only if $\alpha>1/2.$
By analyzing the shape of the image of $\D$ under the mapping $h_\alpha(z)=
k_\alpha(z)/z,$ we obtain the following more refined result.

\begin{thm}\label{thm:geom}
Let $0\le\alpha<1$ and $f\in\K(\alpha).$
Then the following hold:
\begin{enumerate}
\item[(i)]
$
\dfrac{k_\alpha(-r)}{-r}\le\Re\dfrac{f(z)}{z}\le\dfrac{k_\alpha(r)}{r}
$
for $|z|=r<1.$
\item[(ii)]
When $0<\alpha<1/2,$
the asymptotic lines of the boundary curve of $h_\alpha(\D)$ are given by
$v=\pm\cot(\pi\alpha)(u-\frac{1}{2\alpha-1}).$
In particular, the values of $f(z)/z$ for $z\in\D$ are contained in the sector
$S=\{u+iv: |v|< \cot(\pi\alpha)(u-\frac1{2\alpha-1})\}.$
\item[(iii)]
When $1/2\le\alpha<1,$
$$
\left|\Im\frac{f(z)}{z}\right|< M(\alpha),\quad z\in\D,
$$
where
$$
M(\alpha)=\max_{0<\theta<\pi}\Im[e^{-i\theta}k_\alpha(e^{i\theta})]
\le M(\tfrac12)=\frac\pi2.
$$
The estimate is sharp.
\end{enumerate}
\end{thm}

We remark that the left-hand inequality in (i) was already proved
by Brickman {\it et al.} \cite[Theorem 10]{BHMW73} and the right-hand one
follows also from Robertson's theorem (see Lemma \ref{lem:sub} below).
A much simpler proof of (i) is now available thanks to Corollary \ref{cor:main}.
The proof of this theorem and
more information about the constant $M(\alpha)$ will be given in Section 3.
We also provide an application of our results to an extremal problem
for $\K(\alpha)$ in Section 3.

Styer and Wright \cite{SW73} studied (non-)univalence of a convex combination
of two convex functions.
Among other things, the following result is most relevant to the present study.

\begin{Thm}[Styer and Wright]
Let $f, g\in\K$ be odd convex functions.
If $|\Im[f(z)/z]|<\pi/4$ and $|\Im[g(z)/z]|<\pi/4$ on $|z|<1,$
then $(f+g)/2\in\es^*.$
\end{Thm}

Styer and Wright suspected that the assumption $|\Im[f(z)/z]|<\pi/4$
in the theorem was superfluous.
They even stated the belief that
\begin{equation}\label{eq:H2}
\frac{f(z)}z\prec H_2(z)
:=\frac1{2z}\log\frac{1+z}{1-z}
=\sum_{n=0}^\infty\frac{z^{2n}}{2n+1}
\end{equation}
if $f\in\K$ is odd; namely, $f(-z)=-f(z).$
Note that $|\Im H_2(z)|<\pi/4$ on $|z|<1.$
Indeed, Hallenbeck and Ruscheweyh \cite{HR75} proved that
\begin{equation}\label{eq:H1}
\frac{f(z)}{z}\prec H_1(z):=\frac{1}{2\sqrt{z}}\log\frac{1+\sqrt{z}}{1-\sqrt{z}}
=\sum_{n=0}^\infty\frac{z^{n}}{2n+1}
\end{equation}
for a function $f\in\K$ with $f''(0)=0,$
which implies that $|\Im[f(z)/z]|<\pi/4.$
In this way, they strengthened the above theorem
(see \cite[Corollary 2]{HR75}):

\begin{Thm}[Hallenbeck and Ruscheweyh]
Let $f,g\in\K$ satisfy $f''(0)=g''(0)=0.$
Then $(f+g)/2\in\es^*.$
\end{Thm}

We give another result of this type.

\begin{thm}\label{thm:starlike}
$(f+g)/2\in\es^*$ for $f,g\in\K(0.6).$
\end{thm}

The proof will be given in Section 4.
Note that the constant $0.6=3/5$ is not best possible.

We remark that the claim \eqref{eq:H2} for an odd convex function $f$
is not necessarily true.
An example will be given in Section 5.

\section{Proof of Theorem \ref{thm:convex}}

We now show that the function $h_\alpha(z)=k_\alpha(z)/z$
is convex (univalent) on $\D$ for each $0\le\alpha<1.$
To this end, we only need to see that $1+zh_\alpha''(z)/h_\alpha'(z)$
has positive real part.
Since the case $\alpha=0$ is trivial, we assume that $\alpha>0.$
Put $\beta=2-2\alpha\in(0,2)$ for convenience.
We assume $\alpha\ne1/2$ so that $\beta\ne1$ for a while.
A simple calculation yields
$$
h_\alpha'(z)
=\frac{(1-\beta z)(1-z)^{-\beta}-1}{(1-\beta)z^2}
$$
and
\begin{equation}\label{eq:h}
1+\frac{zh_\alpha''(z)}{h_\alpha'(z)}
=-1-\frac{\beta(1-\beta)z^2}{((1-z)^{\beta}-1+\beta z)(1-z)}.
\end{equation}
With the Pochhammer symbol $(a)_n=a(a+1)\cdots(a+n-1),$
we compute
\begin{align*}
(1-z)^{\beta}-1+\beta z
&=\sum_{n=2}^\infty\frac{(-\beta)_n}{(1)_n}z^n \\
&=\frac{-\beta(1-\beta)z^2}{2}\sum_{n=2}^\infty
\frac{(2-\beta)_{n-2}}{(3)_{n-2}}z^{n-2} \\
&=-\frac{\beta(1-\beta)z^2}{2}\sum_{n=0}^\infty
\frac{(2-\beta)_{n}}{(3)_{n}}z^{n}.
\end{align*}
Letting $b_n=(2-\beta)_n/(3)_n$ for $n\ge0,$ we obtain
\begin{align*}
-\frac{((1-z)^{\beta}-1+\beta z)(1-z)}{\beta(1-\beta)z^2}
&=~\frac{1-z}2\sum_{n=0}^\infty b_nz^n \\
&=\frac12\left(1+\sum_{n=1}^\infty (b_n-b_{n-1})z^n\right) \\
&=\frac{1+\omega(z)}2,
\end{align*}
where
$$
\omega(z)=\sum_{n=1}^\infty (b_n-b_{n-1})z^n.
$$
Hence, we have the expression
$$
1+\frac{zh_\alpha''(z)}{h_\alpha'(z)}
=-1+\frac2{1+\omega(z)}
=\frac{1-\omega(z)}{1+\omega(z)}.
$$
Note that this is valid also for $\alpha=1/2$ as is confirmed directly
or by taking limit as $\alpha\to1/2.$

In order to show $\Re(1+zh_\alpha''(z)/h_\alpha'(z))>0,$ it suffices to check
$|\omega(z)|<1.$
Since
$$
\frac{b_{n}}{b_{n-1}}=\frac{n+1-\beta}{n+2}=1-\frac{1+\beta}{n+2}<1,
$$
we see that $\{b_n\}$ is a decreasing sequence of positive numbers.
Therefore,
$$
|\omega(z)|\le\sum_{n=1}^\infty(b_{n-1}-b_n)|z|^n
<\sum_{n=1}^\infty(b_{n-1}-b_n)=b_0-\lim_{n\to\infty}b_n\le b_0=1
$$
for $z\in\D$ as required.
(Indeed, we can easily show that $b_n\to0$ as $n\to\infty.$)
\qed

\bigskip

We remark that the function $k_\alpha$ can be expressed in terms of
the Gauss hypergeometric function
$$
\null_2F_1(a,b;c;z)=\sum_{n=0}^\infty\frac{(a)_n(b)_n}{(c)_n}\cdot\frac{z^n}{n!}.
$$
Indeed, by integrating both sides of
$$
k_\alpha'(z)=(1-z)^{-\beta}=\sum_{n=0}^\infty(\beta)_n\frac{z^n}{n!}
$$
with $\beta=2-2\alpha,$ we obtain
$$
k_\alpha(z)=\sum_{n=0}^\infty\frac{(\beta)_n}{n+1}\cdot\frac{z^{n+1}}{n!}
=z\sum_{n=0}^\infty\frac{(\beta)_n(1)_n}{(2)_n}\cdot\frac{z^n}{n!},
$$
and hence
$$
h_\alpha(z)=\frac{k_\alpha(z)}{z}=\null_2F_1(\beta,1;2;z).
$$
We extract the following result from
K\"ustner's theorems in \cite{Kus02} (Theorem 1.1 with $r=1$ and Remark 2.3,
see also Corollary 6 (a) in \cite{Kus07}).

\begin{lem}[K\"ustner]
For non-zero real numbers $a,b,c$ with $-1<a\le b<c,$
let $F(z)=\null_2F_1(a,b;c;z).$
Then
$$
\inf_{z\in\D}\left(1+\frac{zF''(z)}{F'(z)}\right)
=1+\frac{-F''(-1)}{F'(-1)}
\ge 1-\frac{(a+1)(b+1)}{b+c+2}
$$
\end{lem}

Since $\null_2F_1(a,b;c;z)=\null_2F_1(b,a;c;z),$ we can apply the above
lemma to our function $h_\alpha(z)=\null_2F_1(\beta,1;2;z)$
for $0<\alpha<1;$ equivalently, for $0<\beta<2.$
Hence, by \eqref{eq:h}, we obtain
\begin{align*}
\inf_{z\in\D}\left(1+\frac{zh_\alpha''(z)}{h_\alpha'(z)}\right)
&=1-\frac{h_\alpha''(-1)}{h_\alpha'(-1)}
=\frac{2^{\beta+1}-2-\beta-\beta^2}{2(1+\beta-2^\beta)} \\
&\ge\begin{cases}
\dfrac{4\alpha-1}5, & 1/2\le \alpha<1, \\
\null & \null \\
\dfrac{\alpha}{3-\alpha}, & 0<\alpha\le 1/2.
\end{cases}
\end{align*}
In this way, we have obtained another proof of convexity of $h_\alpha.$

\section{Mapping properties of functions in $\K(\alpha)$}

The present section is devoted to the proof of Theorem \ref{thm:geom}.
Before the proof, we note basic results due to Robertson \cite{Rob36}
(see also Pinchuk \cite{Pin68}).

\begin{lem}[Robertson]\label{lem:sub}
Let $0\leq \alpha<1$ and $f\in\K(\alpha).$
Then,
$$
-k_\alpha(-r)\le |f(z)|\le k_\alpha(r)\quad\text{for}~ |z|=r<1.
$$
In particular,
the image domain $f(\D)$ contains the disk $|w|<-k_\alpha(-1).$
\end{lem}

We will use also the following simple fact.

\begin{lem}\label{lem:conv}
Let $\Omega$ be an unbounded convex domain in $\C$ whose
boundary is parametrized positively by a Jordan curve $w(t)=u(t)+iv(t),~
0<t<1,$ with $w(0^+)=w(1^-)=\infty.$
Suppose that $u(0^+)=+\infty$ and that $v(t)$ has a finite limit
as $t\to0^+.$
Then $v(t)\le v(0^+)$ for $0<t<1.$
\end{lem}

\begin{pf}
Let $0\le t^*\le1$ be the number such that
$u(t^*)=\inf_{0<t<1}u(t)$ and that $u(t)>u(t^*)$ for $0<t<t^*.$
(We interpret $u(0)=u(0^+)$ or $u(1)=u(1^-)$
when $t^*=0$ or $1,$ respectively.)
By the assumption $u(0^+)=+\infty,$ we have $t^*>0.$
Note that $u(t)$ is strictly decreasing in $0<t<t^*.$
By convexity and orientation, the part $w((t^*,1))$ of the boundary lies below
the part $w((0,t^*)).$
Thus, it is enough to show that $v(t)$ is non-increasing in $0<t<t^*.$
Let $0<t_0<t_1<t_2<t^*$ and set $w(t_j)=u_j+iv_j$ for $j=0,1,2.$
By convexity, the part $w((t_0,t_2))$ of the boundary lies above
the line which passes through the points $w(t_0)$ and $w(t_2);$
equivalently,
$$
v(t)\ge v_2+\frac{v_0-v_2}{u_0-u_2}(u(t)-u_2),\quad t_0<t<t_2.
$$
We now put $t=t_1$ and let $t_0\to0^+$ to obtain $v_1=v(t_1)\ge v_2=v(t_2).$
Thus we have shown that $v(t)$ is non-increasing as required.
\end{pf}

We are now ready to prove Theorem \ref{thm:geom}.

\begin{pf}[Proof of Theorem \ref{thm:geom}]
Since $h_\alpha(z)=k_\alpha(z)/z$ is convex and symmetric in $\R,$
we easily see that $h_\alpha(-r)\le \Re h_\alpha(z)\le h_\alpha(r)$
for $|z|=r<1.$
Therefore, assertion (i) immediately follows from Corollary \ref{cor:main}.

To prove (ii) and (iii), we study mapping properties of the function 
$h_\alpha(z).$
We remark that $h_\alpha$ analytically extends to $\partial\D\setminus\{1\}$
by its form.
Let us investigate the shape of the boundary of $h_\alpha(\D).$
In the rest of this section, it is convenient to put $\gamma=2\alpha-1\in[-1,1).$
Note that $\gamma<0$ if and only if $\alpha<1/2.$
We write $h_\alpha(e^{i\theta})=u_\gamma(\theta)+iv_\gamma(\theta)$
for $0<\theta<2\pi.$
We remark that the symmetry $h_\alpha(\bar z)=\overline{h_\alpha(z)}$
leads to the relations $u_\gamma(2\pi-\theta)=u_\gamma(\theta)$
and $v_\gamma(2\pi-\theta)=-v_\gamma(\theta).$
Thus, we may restrict our attention to the range $0<\theta\le\pi.$
It is easy to obtain the following expressions for $\gamma\ne0$:
\begin{align*}
u_\gamma(\theta)&=\frac{-1}\gamma
\left(\left(2\sin\frac{\theta}{2}\right)^{\gamma}
\cos\left(-\theta+\frac{\theta-\pi}{2}\gamma\right)-\cos\theta\right), \\
v_\gamma(\theta)&=\frac{-1}\gamma
\left(\left(2\sin\frac{\theta}{2}\right)^{\gamma}
\sin\left(-\theta+\frac{\theta-\pi}{2}\gamma\right)+\sin\theta\right).
\end{align*}
Observe that for $-1<\gamma<0,$ both $u_\gamma(\theta)$ and
$v_\gamma(\theta)$ tend to $+\infty$ as $\theta\to 0^+.$
A simple calculation yields
$$
\lim_{\theta\to0^+}\frac{v_\gamma(\theta)}{u_\gamma(\theta)}
=\tan\frac{-\pi\gamma}{2}
=-\tan\frac{\pi\gamma}{2}
$$
and
\begin{align*}
\lim_{\theta\to0^+}
\left(v_\gamma(\theta)+u_\gamma(\theta)\tan\frac{\pi\gamma}2 \right)
&=
\lim_{\theta\to0^+}\left[
\frac{-(2\sin\frac{\theta}{2})^{\gamma}\sin(\frac\gamma2-1)\theta}%
{\gamma\cos(\pi\gamma/2)}-\frac1\gamma\tan\frac{\pi\gamma}2\cos\theta
\right] \\
&=-\frac1\gamma \tan\frac{\pi\gamma}2.
\end{align*}
Therefore, 
\begin{align*}
v&=-\tan\frac{\pi\gamma}2\left(u-\frac1\gamma\right) \\
&=\cot(\alpha\pi)\left(u-\frac{1}{2\alpha-1}\right)
\end{align*}
is an asymptotic line of the boundary curve $\partial h_\alpha(\D).$
Since $h_\alpha(\D)$ is a convex domain symmetric in the real axis,
we conclude assertion (ii).

Next we assume $\alpha\ge1/2$ to show (iii).
Since $f(z)/z\prec k_\alpha(z)/z\prec k_{1/2}(z)/z$ for $f\in\K(\alpha),$
the assertion is clear except for $M(1/2)=\pi/2.$
A simple computation gives us the expression
$$
v_{1/2}(\theta)=\frac{\pi-\theta}2\cos\theta+\sin\theta
\log\left(2\sin\frac\theta2\right)
$$
for $0<\theta<\pi.$
We easily get $v_{1/2}(0^+)=\pi/2.$
Thus we conclude that $M(1/2)=\pi/2$ by Lemma \ref{lem:conv}.
We have thus proved assertion (iii).
\end{pf}


We indicate how to compute the value of $M(\alpha)$ for $1/2<\alpha<1.$
Set $c=\gamma/2=\alpha-1/2\in(0,1/2).$
Since $h_\alpha(\D)$ is a bounded convex domain symmetric in $\R,$
it is easy to see that $v_\gamma(\theta)$ has a unique critical point,
say, $\theta_\alpha$ at which $v_\gamma$ attains its maximum so that
$M(\alpha)=v_\gamma(\theta_\alpha).$
Here, $\theta=\theta_\alpha$ is a unique solution of the equation
\begin{equation}\label{eq:v'}
\left[c\cot\frac\theta2+(1-c)\cot(c\pi+(1-c)\theta)\right]
\left(2\sin\frac{\theta}{2}\right)^{2c}\sin(c\pi+(1-c)\theta)
-\cos\theta=0
\end{equation}
in $0<\theta<\pi,$ where $c=\alpha-1/2.$
By using this equation, we can express $M(\alpha)$ in a different way:
\begin{equation}\label{eq:M}
M(\alpha)=\frac1{2c}\left[
\frac{\cos\theta_\alpha}{c\cot(\theta_\alpha/2)
+(1-c)\cot(c\pi+(1-c)\theta_\alpha)}
-\sin\theta_\alpha\right].
\end{equation}
This expression will be used in the proof of Theorem \ref{thm:starlike}.

Assertion (ii) of Theorem \ref{thm:geom} can
be applied to an extremal problem for $\K(\alpha).$
For $0\le\alpha<1$ and $t\in\R,$ we consider the quantity
$$
Q_\alpha(t)=\inf_{f\in\K(\alpha),\, z\in\D}\Re\left[
e^{it}\frac{f(z)}{z}
\right].
$$
The quantity $M(\alpha)$ in Theorem \ref{thm:geom} is a particular case
of this quantity. Indeed, we have $Q_\alpha(\pi/2)=-M(\alpha)$ for
$1/2\le\alpha<1.$
We have the obvious monotonicity
$Q_\alpha(t)\le Q_\beta(t)$ for $0\le\alpha<\beta<1$
and the symmetry $Q_\alpha(-t)=Q_\alpha(t).$
It is thus enough to consider the case when $0\le t\le \pi.$

\begin{thm}
For $0<\alpha<1,$
the function $\varphi_\alpha(\theta)=\theta+\arg h_\alpha'(e^{i\theta})$ maps
the interval $(0,\pi]$ onto $(\pi(1-\alpha), \pi]$ homeomorphically.
Furthermore, the following hold.
\begin{enumerate}
\item[(i)]
Suppose $\alpha=0.$ Then, $Q_0(0)=1/2$ and $Q_0(t)=-\infty$
for $0<t\le\pi.$
\item[(ii)]
Suppose $0<\alpha<1/2.$ Then
$$
Q_\alpha(t)=\begin{cases}
\Re[e^{i(t-\theta_0)}k_\alpha(e^{i\theta_0})],
& 0\le t<\alpha\pi, \theta_0=\varphi_\alpha\inv(\pi-t), \\
(2\alpha-1)^{-1}\cos(\alpha\pi),
& t=\alpha\pi, \\
-\infty, & \alpha\pi<t\le\pi.
\end{cases}
$$
\item[(iii)]
Suppose $\alpha=1/2.$ Then
$$
Q_{1/2}(t)=\begin{cases}
\Re[e^{i(t-\theta_0)}k_{1/2}(e^{i\theta_0})],
& 0\le t<\pi/2, \theta_0=\varphi_{1/2}\inv(\pi-t), \\
-\pi/2,
& t=\pi/2, \\
-\infty, & \pi/2<t\le\pi.
\end{cases}
$$
\item[(iv)]
Suppose $1/2<\alpha<1.$ Then
$$
Q_\alpha(t)=\begin{cases}
\Re[e^{i(t-\theta_0)}k_\alpha(e^{i\theta_0})],
& 0\le t<\alpha\pi, \theta_0=\varphi_\alpha\inv(\pi-t), \\
(2\alpha-1)\inv\cos t,
& \alpha\pi\le t\le\pi.
\end{cases}
$$
\end{enumerate}
\end{thm}

\begin{pf}
When $t=0$ or $\pi,$ the assertions are clear.
Assume therefore that $0<t<\pi.$
Let $D_\alpha=h_\alpha(\D).$
By Corollary \ref{cor:main}, we have
$$
Q_\alpha(t)=\inf_{u+iv\in D_\alpha}\Re\left[ e^{it}(u+iv)\right]
=\inf_{u+iv\in D_\alpha}\big[u\cos t-v\sin t\big].
$$
Then, geometrically,  we can say that $-Q_\alpha(t)/\sin t$ is the supremum
of $y$-intercepts of those lines $y=x\cot t+C$ which intersect with $D_\alpha.$
Since $D_\alpha$ does not intersect the $y$-axis,
such a line must intersect with $\partial D_\alpha.$
Therefore, in the above characterization of $Q_\alpha(t),$ 
$D_\alpha$ may be replaced by $\partial D_\alpha.$
Hence, noting also the symmetry of $D_\alpha$ in $\R,$ we further obtain
\begin{align*}
Q_\alpha(t)
&=\inf_{u+iv\in\partial D_\alpha}\big[u\cos t-v\sin t\big] \\
&=\inf_{u+iv\in\partial D_\alpha, v\ge0}\big[u\cos t-v\sin t\big] \\
&=\inf_{0<\theta<\pi}F(\theta),
\end{align*}
where
$$
F(\theta)=u_\gamma(\theta)\cos t-v_\gamma(\theta)\sin t.
$$
and $u_\gamma, v_\gamma$ are the functions given by
$h_\alpha(e^{i\theta})=u_\gamma(\theta)+iv_\gamma(\theta)$
with $\gamma=2\alpha-1,$ as before.

When $\alpha=0,$ the function $h_0(z)=1/(1+z)$ maps the unit disk onto
the half-plane $\Re w>1/2$ so that assertion (i) is obvious.
We thus assume that $0<\alpha<1$ in the rest of the proof.

First we analyze the case when $Q_\alpha(t)=-\infty.$
Recall that $u_\gamma(\theta)\to+\infty$ and
$v_\gamma(\theta)=u_\gamma(\theta)\cot(\alpha\pi)+O(1)$ as
$\theta\to0^+$ for $0<\alpha<1/2$ by (ii) of Theorem \ref{thm:geom}.
This is valid also for $\alpha=1/2.$
Hence,
$$
\sin(\alpha\pi)\big[u_\gamma(\theta)\cos t-v_\gamma(\theta)\sin t\big]
=u_\gamma(\theta)\sin(\alpha\pi-t)+O(1)
\to-\infty\quad(\theta\to0^+),
$$
whenever $\sin(\alpha\pi-t)<0,$ which confirms the assertion
for $\alpha\pi<t<\pi$ and $0<\alpha\le1/2.$

We now show the first assertion of the theorem.
Let $\psi_\alpha(\theta)
=\arg[u_\gamma'(\theta)+iv_\gamma'(\theta)]\in(\pi/2,3\pi/2]$
for $0<\theta\le\pi.$
The strict convexity of $D_\alpha$ implies that $\psi_\alpha$ is 
strictly increasing.
Note that 
$
\psi_\alpha(\theta)=\arg h_\alpha'(e^{i\theta})+\theta+\pi/2
=\varphi_\alpha(\theta)+\pi/2.
$
Then we consider the case $0\le t<\alpha\pi.$
From the proof of assertion (ii) of Theorem \ref{thm:geom},
we see that $\psi_\alpha(0^+)=3\pi/2-\alpha\pi$ for $0<\alpha<1/2.$
This is valid also for $1/2\le\alpha<1.$
Indeed, it follows from
$$
\tan\psi_\alpha(0^+)=
\lim_{\theta\to0^+}\frac{v_\gamma(\theta)}{u_\gamma(\theta)-u_\gamma(0^+)}
=-\tan\frac{\pi\gamma}{2}=\cot(\alpha\pi)
$$
for $1/2<\alpha<1.$
We can also see that $\psi_{1/2}(0^+)=\pi$ directly.
Hence, we conclude that the range of $\psi_\alpha(\theta)$ on $0<\theta\le\pi$
is precisely $(\frac{3\pi}2-\alpha\pi,\frac{3\pi}2],$ which proves
the required assertion.

We now consider the case when $0\le t<\alpha\pi.$
Then $F'(\theta)$ vanishes precisely when 
$\tan\psi_\alpha(\theta)=v_\gamma'(\theta)/u_\gamma'(\theta)
=\cot t=\tan(3\pi/2-t),$ namely, $\varphi_\alpha(\theta)=\pi-t.$
Thus we see that $F(\theta)$ takes its minimum at 
$\theta_0=\varphi_\alpha\inv(\pi-t)$ and the corresponding assertions hold.

Our next task is to consider
the borderline case $t=\alpha\pi.$
When $0<\alpha<1/2,$
Theorem \ref{thm:geom} (ii) implies that the supremum of the $y$-intercepts
of the lines $y=x\cot(\alpha\pi)+k$ intersecting with $D_\alpha$ is
$\cot(\alpha\pi)/(1-2\alpha).$ This case has been confirmed to be true.
When $\alpha=1/2,$ the assertion is contained in Theorem \ref{thm:geom} (iii).
When $\alpha>1/2,$ this case can be included in the final case below.

We finally consider the case when $1/2<\alpha<1$ and $\alpha\pi\le t<\pi.$
In this case the function $F(\theta)$ has no critical point in 
$0<\theta<\pi.$
Since $F'(\pi)=-v_\gamma'(\pi)\sin t>0,$ we see that $F(\theta)$ is increasing
in $0<\theta<\pi$ so that $Q_\alpha(t)=F(0^+)=(2\alpha-1)\inv\cos t.$
\end{pf}

\section{Proof of Theorem \ref{thm:starlike}}

We denote by $\D_r$ the disk $|z|<r.$
Throughout this section, we define $f_a$ for $f\in\A_1$ and $a\in\D$ by
$f_a(z)=f(az)/a.$
Here, we set $f_0(z)=\lim_{a\to0}f_a(z)=z.$
We begin with the following simple observation.

\begin{lem}\label{lem:covering}
Let $f\in\es.$
Suppose that $f(\D)$ contains the disk $\D_\rho$ for some $\rho>0.$
Then $\D_\rho\subset f_a(\D)$ for $a\in\D.$
\end{lem}

\begin{pf}
It suffices to show that $\D_{\rho r}\subset f(\D_r)$ for $0<r<1.$
By assumption, $g(w)=f\inv(\rho w)$ is a univalent analytic function
on $\D$ with $|g(w)|<1$ and $g(0)=0.$
Then the Schwarz lemma implies that $g(\D_r)\subset\D_r,$
which in turn gives us $\D_{\rho r}\subset f(\D_r)$ as required.
\end{pf}

By making use of the idea due to Styer and Wright \cite{SW73},
the following result can now be shown.
For convenience of the reader, we reproduce the proof here in a
somewhat simplified form.

\begin{lem}
Let $\rho$ be a positive constant.
Suppose that two functions $f$, $g\in\K$ satisfy the following two conditions:
\begin{enumerate}
\item
$f(\D)$ and $g(\D)$ both contain the disk $\D_\rho,$ and
\item
$|\Im[f(z)/z]|<\rho$ and $|\Im[g(z)/z]|<\rho$ on $\D.$
\end{enumerate}
Then $(f+g)/2\in\es^*.$
\end{lem}

\begin{pf}
Put $h=f+g.$
For starlikeness, we need to show that $\Re[zh'(z)/h(z)]>0$ on $\D.$
We will show that $\Re[zf'(z)/h(z)]>0.$
Since we can do the same for $g,$ it will finish the proof.

Let $a\in\D$ with $a\ne0.$
Since $f_a'(1)/h_a(1)=af'(a)/h(a),$ it is enough to show the inequality
$\Re[f_a'(1)/h_a(1)]\ge0.$
Denote by $W$ the set $\{w: |w|\ge\rho, |\Im w|<\rho\}.$
Then $W$ consists of the two connected components $W_+$ and $W_1,$
where $W_\pm=\{w\in W: \pm \Re w>0\}.$
By Lemma \ref{lem:covering} and the relation $f_a(z)/z=f(az)/(az),$
the assumptions imply $f_a(1)\in W.$
Since the (continuous) curve $t\mapsto f_{ta}(1),~0\le t\le 1,$ connects $f_a(1)$
with $f_0(1)=1,$ we see that $f_a(1)\in W_+.$
Since we have $g_a(1)\in W_+$ in the same way and thus $-g_a(1)\in W_-,$
the segment $[-g_a(1), f_a(1)]$ intersects the disk $\D_\rho.$
Choose a point $w_0\in[-g_a(1), f_a(1)]\cap\D_\rho.$
Then the vector $f_a(1)-w_0$ is directed at the point $f_a(1)$
outward from the convex domain $f_a(\D).$
Since the tangent vector of the curve $f_a(e^{i\theta})$ at $\theta=0$ is 
given by $if_a'(1),$ we have
$$
\arg[if_a'(1)]-\pi\le \arg[f_a(1)-w_0]=\arg[f_a(1)+g_a(1)]\le
\arg[if_a'(1)],
$$
which is equivalent to $|\arg[f_a'(1)/h_a(1)]|\le\pi/2.$
Thus we have shown the desired inequality $\Re[f_a'(1)/h_a(1)]\ge0.$
\end{pf}

\begin{pf}[Proof of Theorem \ref{thm:starlike}]
Let $f,g\in\K(3/5).$
We will apply the last lemma to these two functions.
Let $\rho= -k_{3/5}(-1)=5(2^{1/5}-1)=0.743491\dots.$
By Theorem \ref{thm:geom}, we have only to show that $M(3/5)\le\rho.$
We denote by $F(\theta)$ the function in the left-hand side in \eqref{eq:v'}
for $c=\frac35-\frac12=\frac1{10}.$
A numerical computation gives us $F(0.11)=0.0050\dots>0$ and
$F(0.114)=-0.0010\dots<0.$
Thus we have $0.11<\theta_{3/5}<0.114.$
By \eqref{eq:M}, we have the expression
$M(3/5)=5G(\theta_{3/5}),$ where
$$
G(\theta)=\frac{\cos\theta}{c\cot(\theta/2)
+(1-c)\cot(c\pi+(1-c)\theta)}-\sin\theta
=\frac{\cos\theta}{H(\theta)}-\sin\theta.
$$
We observe that $H(\theta)$ is positive and decreasing in
$0<\theta<\frac{1/2-c}{1-c}\pi=4\pi/9,$
because
$$
H'(\theta)=-\frac{c}{2\sin^2(\theta/2)}-\frac{(1-c)^2}{\sin^2(c\pi+(1-c)\theta)}
<0.
$$
Also, we see that $-H'(\theta)$ is positive and decreasing in 
$0<\theta<4\pi/9$ by its form.
Since
$$
G'(\theta)=-\frac{\sin\theta}{H(\theta)}
-\frac{H'(\theta)\cos\theta}{H(\theta)^2}-\cos\theta,
$$
letting $\theta_0=0.11$ and $\theta_1=0.114,$
we estimate on $\theta_0\le\theta\le\theta_1$ in the form
$$
G'(\theta)>-\frac{\sin\theta_1}{H(\theta_1)}
-\frac{H'(\theta_1)\cos\theta_1}{H(\theta_0)^2}-\cos\theta_0
=0.326\dots>0.
$$
Hence $G(\theta)$ is increasing in this interval so that
$$
M(3/5)=5G(\theta_{3/5})<5G(\theta_1)=0.743487\dots<\rho.
$$
The proof is now complete.
\end{pf}

\section{An example}

We conclude the present note by giving an example of an odd convex
function $f\in\K$ such that 
$$
\frac{f(z)}{z}\not\prec
H_2(z)=\frac{1}{2z}\log\frac{1+z}{1-z}\quad\text{on}~\D.
$$

The following result due to Alexander \cite{Alexander15}
(see also Goodman \cite{Good57}) is useful for our aim here.

\begin{lem}[Alexander]
The function $f(z)=z+a_2z^2+a_3z^3+\dots$ is convex univalent on $\D$ if
$$
\sum_{n=1}^\infty n^2|a_n|\le 1.
$$
\end{lem}

We also need the following auxiliary result which is a special case
of Theorem 5 of Ruscheweyh \cite{Rus75} with $n=1.$

\begin{lem}[Ruscheweyh]
The function
$$
q_\gamma(z)=\sum_{j=1}^\infty\frac{\gamma+1}{\gamma+j}z^j
$$
belongs to $\K$ for $\Re\gamma\ge0.$
\end{lem}

In particular, the function $H_1$ given in \eqref{eq:H1} is
univalent because $H_1=1+q_{1/2}/3.$

We now consider the function
$$
f(z)=z+\frac{z^3}{100}+\frac{z^5}{50}.
$$
Then, by Alexander's lemma, $f$ is an odd convex function.
Secondly, we observe that $f$ has a non-zero fixed point $z_0$ in $\D.$
Indeed, by solving the algebraic equation $f(z)=z,$ we obtain
$z_0=\pm i/\sqrt2.$

We now show that $f(z)/z$ is not subordinate to $H_2(z)$
given in \eqref{eq:H2}.
Suppose, to the contrary, that
$$
\frac{f(z)}{z}\prec
H_2(z)=\frac{1}{2z}\log\frac{1+z}{1-z},\quad z\in\D.
$$
Then there exists an analytic function $\omega$ on $\D$ with $\omega(0)=0$
and $|\omega|<1$ such that
$$
\frac{f(z)}{z}=H_2(\omega(z))=H_1(\omega(z)^2).
$$
Thus
$$
\frac{zf'(z)-f(z)}{z^2}=2\omega(z)\omega'(z)H_1'(\omega(z)^2).
$$
Since $f(z_0)=z_0,$ we have $H_1(\omega(z_0)^2)=1=H_1(0).$
Univalence of $H_1$ enforces the relation $\omega(z_0)=0$ to hold.
Hence, $z_0f'(z_0)-f(z_0)=0$ which is equivalent to $f'(z_0)=1.$
By solving the equation $f'(z)=1,$ we obtain
$z_0=\pm i\sqrt{3/10}.$
This is a contradiction.
Therefore, $f(z)/z$ is not subordinate to $H_2(z).$


\def\cprime{$'$} \def\cprime{$'$} \def\cprime{$'$}
\providecommand{\bysame}{\leavevmode\hbox to3em{\hrulefill}\thinspace}
\providecommand{\MR}{\relax\ifhmode\unskip\space\fi MR }
\providecommand{\MRhref}[2]{%
  \href{http://www.ams.org/mathscinet-getitem?mr=#1}{#2}
}
\providecommand{\href}[2]{#2}

\end{document}